\documentclass[11pt,reqno]{amsart}

\usepackage[T1]{fontenc}
\usepackage[utf8]{inputenc}
\usepackage{lmodern}
\usepackage{amsmath,amssymb,amsthm,mathtools}
\usepackage{enumitem}
\usepackage{microtype}
\usepackage{needspace}
\usepackage{xcolor}
\usepackage[
  colorlinks=true,
  linkcolor=blue!45!black,
  citecolor=blue!45!black,
  urlcolor=blue!45!black
]{hyperref}

\hypersetup{
  pdftitle={An Explicit Uniform Connectivity Bound for Generalized Markoff Graphs},
  pdfauthor={Mingxu Jiang}
}

\allowdisplaybreaks[2]
\numberwithin{equation}{section}

\newtheorem{theorem}{Theorem}[section]
\newtheorem{proposition}[theorem]{Proposition}
\newtheorem{lemma}[theorem]{Lemma}

\title[Uniform connectivity of generalized Markoff graphs]
{An Explicit Uniform Connectivity Bound for Generalized Markoff Graphs}

\author{Mingxu Jiang}
\address{Department of Mathematical and Statistical Sciences,
Clemson University, Clemson, SC, USA}
\email{mingxuj@clemson.edu}


\begin{document}

\begin{abstract}
Let $X_D$ be the surface $x^2+y^2+z^2-xyz=D$ over $\mathbb F_p$.
For every prime $p>10^{320}$ and every
$D\in\mathbb F_p\setminus\{4\}$, we prove that the three Vieta
involutions act transitively on the trace triples of generating pairs
of $\mathrm{SL}_2(\mathbb F_p)$.
The proof combines an explicit primitive-fiber argument, Chen's
orbit-size lower bound, and a weighted maximal-divisor estimate.
We also give an explicit descent from Nielsen transitivity to
transitivity under the Vieta involutions alone.
\end{abstract}

\maketitle

\section{Introduction}\label{sec:introduction}

Let $p$ be an odd prime and let $D\in\mathbb F_p$.
We consider the generalized Markoff surface
\[
X_D:\quad x^2+y^2+z^2-xyz=D
\]
and the three Vieta involutions
\[
\begin{aligned}
s_x(x,y,z)&=(yz-x,y,z),\\
s_y(x,y,z)&=(x,xz-y,z),\\
s_z(x,y,z)&=(x,y,xy-z).
\end{aligned}
\]
Set $\Gamma=\langle s_x,s_y,s_z\rangle$ and
$H=\langle\Gamma,S_3\rangle$, where $S_3$ permutes coordinates.
For $D\ne4$, define the \emph{essential locus}
\[
S_D(p)=
\left\{
(\operatorname{tr}A,\operatorname{tr}B,\operatorname{tr}(AB)):
\begin{array}{l}
A,B\in\mathrm{SL}_2(\mathbb F_p),\
\langle A,B\rangle=\mathrm{SL}_2(\mathbb F_p),\\
\operatorname{tr}[A,B]=D-2
\end{array}
\right\}.
\]
The generalized Markoff graph considered here has vertex set $S_D(p)$
and edges joining each point to its images under $s_x,s_y,s_z$.
Its connectivity is equivalent to $\Gamma$-transitivity.

\Needspace{5\baselineskip}
\begin{theorem}[Uniform pure Vieta transitivity]\label{thm:main}
For every prime $p>10^{320}$ and every
$D\in\mathbb F_p\setminus\{4\}$, the pure Vieta group $\Gamma$
acts transitively on $S_D(p)$.
\end{theorem}

The classical case $D=0$ belongs to the strong approximation
conjecture posed by Baragar \cite{Baragar91}; its generalization to
essential trace triples is the $Q$-classification conjecture of
McCullough and Wanderley \cite{MW13}.
For $D=0$ and $p\ge7$, their classification of nonessential triples
gives
$S_0(p)=X_0(\mathbb F_p)\setminus\{(0,0,0)\}$
\cite[\S~11]{MW13}.
Bourgain, Gamburd, and Sarnak developed the cage and endgame
framework using rotations, proved that a single orbit has complement of size
at most $p^\varepsilon$ for every fixed $\varepsilon>0$ and all
sufficiently large $p$, and established transitivity outside a sparse
set of primes \cite{BGS26}.
Chen then proved the divisibility of every orbit size by $p$,
which, combined with the Bourgain--Gamburd--Sarnak theorem, gives
transitivity for all sufficiently large primes
\cite[Theorem~1.2.5]{Chen24}.
Martin subsequently gave an elementary proof of Chen's divisibility
theorem \cite{MartinChen25}.
Eddy, Fuchs, Litman, Martin, and Tripeny made the connectivity
threshold explicit and introduced maximal divisors to obtain
$p>3.45\cdot10^{392}$ \cite[Theorem~1.4]{EFLMT25}.

For general $D\ne4$, Bourgain, Gamburd, and Sarnak also announced
extensions of their strong approximation methods to Markoff-type
surfaces \cite{BGS26}.
Kingsbury-Neuschotz developed this framework for the relative
character varieties of four-times punctured spheres, including the
surfaces $X_D$, and proved uniform small-complement estimates
\cite{KN26}.
Chen's component-degree theorem combines with these estimates to
give transitivity for all sufficiently large primes, uniformly in $D$
\cite{Chen26}.
Independently, Martin proved pure Vieta transitivity when
$24\,504\,480\nmid p^2-1$ \cite[Theorem~1.1]{Martin25}.
Our purpose is to make the resulting uniform threshold completely
explicit, using the geometric strategy of Bourgain, Gamburd, and
Sarnak and Chen's lower bound for orbit sizes.

Our arithmetic input refines the maximal-divisor method of
Eddy et al.\ \cite{EFLMT25}, treating
the prime factors $2,3$ jointly and exploiting the coprime odd parts
of $p-1$ and $p+1$.
Following their terminology, a divisor $d$ of $n$ is
\emph{maximal with respect to $m$} if $d\le m$ and no strictly larger
multiple of $d$ dividing $n$ is at most $m$.
Thus, for an integer $n\ge1$ and real $m\ge1$, the set of maximal
divisors is
\[
\mathcal M_m(n)=
\bigl\{
d\mid n:d\le m,\quad
\nexists e\mid n\ \text{with }d\mid e
\text{ and }d<e\le m
\bigr\}.
\]
Set
\[
K_p(m)=|\mathcal M_m(p-1)|+|\mathcal M_m(p+1)|.
\]
The following weighted estimate is the arithmetic refinement used
in the proof of Theorem~\ref{thm:main}.

\Needspace{12\baselineskip}
\begin{theorem}[Weighted maximal-divisor estimate]
\label{thm:antichain}
Set
\[
s=\frac9{100},\qquad
t=\frac{67}{1000},\qquad
K_0=19\,400\,000.
\]
For every odd integer $P\ge3$, every real $m\ge1$, and every $\eta\ge0$,
\begin{equation}\label{eq:weighted-antichain}
\sum_{\epsilon=\pm1}
\ \sum_{d\in\mathcal M_m(P+\epsilon)}d^\eta
<K_0(P+1)^t m^{\eta+s}.
\end{equation}
In particular, for each odd prime $p$,
\begin{equation}\label{eq:K-bound}
K_p(m)<K_0(p+1)^t m^s.
\end{equation}
\end{theorem}

Let $\chi_p$ be the quadratic character of $\mathbb F_p$.
Throughout the geometric and numerical arguments we select
\begin{equation}\label{eq:choose-q}
q=p-\chi_p(-1)\in\{p-1,p+1\},\qquad 4\mid q.
\end{equation}
Write $\omega(q)$ for the number of distinct prime factors of $q$
and $\varphi(q)$ for Euler's totient function, and set
\begin{equation}\label{eq:JB}
J(q)=\frac{2^{\omega(q)}q}{\varphi(q)},\qquad
B_p=\frac{p+1}{p}J(q)(4\sqrt p+8).
\end{equation}
A \emph{good primitive trace} is a value
$a=\lambda+\lambda^{-1}\in\mathbb F_p$ with
$\lambda\in\mathbb F_{p^2}^{\times}$,
$\operatorname{ord}(\lambda)=q$, and $a^2\ne D$.
For $j\in\{1,2,3\}$, define the \emph{coordinate fibers} and their
primitive union by
\[
\begin{aligned}
C_j(a)&=\{(x_1,x_2,x_3)\in X_D(\mathbb F_p):x_j=a\},\\
\mathcal C&=
\bigcup_{a\ \text{good primitive}}\ \bigcup_{j=1}^3 C_j(a).
\end{aligned}
\]
Thus $\mathcal C$ consists of the points having at least one good
primitive coordinate.
For $D=0$, this is the part of the Bourgain--Gamburd--Sarnak cage
construction obtained by choosing just the semisimple torus of
order $q$, rather than both torus types and the available parabolic
fibers \cite[\S~3.2]{BGS26}.
The exclusion $a^2=D$ removes degenerate fibers; the rotation
parameterization in Section~\ref{sec:prelim} will show that a good
primitive fiber is a single rotation orbit.
Proposition~\ref{prop:endgame} will connect these fibers in one
$H$-orbit $\mathcal U\subseteq S_D(p)$ under the hypothesis below.

Let $V\cong C_2^2$ be the group of even coordinate sign changes,
\[
V=\{(x,y,z)\mapsto(\epsilon_1x,\epsilon_2y,\epsilon_3z):
\epsilon_i\in\{\pm1\},\ \epsilon_1\epsilon_2\epsilon_3=1\}.
\]
The choice $4\mid q$ makes $\mathcal U$ invariant under $V$.
Chen's lower bound on Nielsen orbit sizes
among sign classes \cite[Theorem~1.1]{Chen26},
together with the freeness of the \(V\)-action on
\(S_D(p)\), gives the following lower bound for
the complement of \(\mathcal U\).

\Needspace{12\baselineskip}
\begin{theorem}[Lower bound for the complement]\label{thm:four-p}
Suppose that $p\ge113$, $D\ne4$, and
\begin{equation}\label{eq:endgame-condition}
p>J(q)(16\sqrt p+16).
\end{equation}
Then $\mathcal U=H\mathcal C\subseteq S_D(p)$ is a single $H$-orbit
invariant under $V$.
If $\mathcal B=S_D(p)\setminus\mathcal U$ is nonempty, then
\[
|\mathcal B|\ge4p.
\]
The hypotheses hold for every prime $p>10^{12}$ and every $D\ne4$.
\end{theorem}

Finally, the Nielsen action on trace triples allows coordinate
permutations, whereas Theorem~\ref{thm:main} uses only Vieta
involutions.
The following descent makes explicit Martin's normality and
fixed-point argument \cite[proof of Theorem~4.7]{Martin25}, together
with the odd-order rotation argument of Kingsbury-Neuschotz
\cite[Lemmas~7.7--7.8]{KN26}.

\Needspace{5\baselineskip}
\begin{theorem}[Descent to pure Vieta involutions]\label{thm:descent}
Let $p\ge113$ and $D\ne4$.
If $H=\langle\Gamma,S_3\rangle$ acts transitively on $S_D(p)$,
then so does $\Gamma$.
\end{theorem}

The proof is divided into geometric and arithmetic parts.
Section~\ref{sec:prelim} recalls the trace correspondence and rotations.
Section~\ref{sec:geometry} establishes the primitive-fiber endgame,
sign invariance, and descent.
Section~\ref{sec:arithmetic} supplies the torsion-intersection estimate,
the maximum-order obstruction, and Theorem~\ref{thm:antichain}.
Section~\ref{sec:proof} combines the estimates.
Appendix~\ref{app:certificate} records the finite numerical verifications.

\section{Trace representations and rotations}\label{sec:prelim}

The Nielsen action of $\operatorname{Out}(F_2)$ on trace triples has
image $H$. This description appears in McCullough \cite[\S~1]{McC05} and
McCullough and Wanderley
\cite[\S~8]{MW13}; see also Martin \cite[(1.7)]{Martin25}.
The orientation-preserving subgroup is
$\operatorname{Out}^+(F_2)\cong\mathrm{SL}_2(\mathbb Z)$.
In particular, $H$ preserves $S_D(p)$.
A point with two zero coordinates is called \emph{axial}; a point
of $X_D(\mathbb F_p)$ outside $S_D(p)$ is called \emph{nonessential}.

\begin{proposition}[Literature inputs]\label{prop:literature}
Let $p\ge5$ and $D\ne4$.
\begin{enumerate}[label=(\arabic*)]
\item Every point of $X_D(\mathbb F_p)$ is the trace triple of an
absolutely irreducible pair in $\mathrm{SL}_2(\mathbb F_p)$,
unique up to simultaneous $\mathrm{GL}_2(\mathbb F_p)$-conjugacy.

\item Axial triples are nonessential.
Every nonessential nonaxial triple has all three coordinates among
the roots of
\begin{equation}\label{eq:exceptional-poly}
T(T^2-1)(T^2-2)(T^4-3T^2+1),
\end{equation}
and its coordinate eigenvalue orders belong to
\begin{equation}\label{eq:exceptional-orders}
\{3,4,5,6,8,10\}.
\end{equation}

\item The group $V$ preserves $S_D(p)$, and the Nielsen action of
$\operatorname{Out}(F_2)$ descends to $S_D(p)/V$. Every
$\operatorname{Out}^+(F_2)$-orbit $\Omega$ in $S_D(p)/V$
satisfies $|\Omega|\ge p$.
\end{enumerate}
\end{proposition}

Part (1) goes back to Macbeath \cite[Theorems~1--3]{Macbeath69};
see McCullough and Wanderley \cite[\S\S~7--8]{MW13} and Chen
\cite[Theorem~5.2.10]{Chen24} for the formulation used here.
Part (2) follows from Dickson's subgroup classification
\cite{Dickson01} and the trace analysis of Macbeath \cite{Macbeath69}.
The exceptional trace criteria go back to Churchill over $\mathbb C$
\cite{Churchill99} and McCullough over general fields \cite{McC05};
they enter McCullough and Wanderley's classification
\cite[\S~11]{MW13}.
Later formulations include Martin \cite[Theorem~2.5]{Martin25},
Campos-Vargas \cite[Theorem~6]{CV25}, and Chen
\cite[Lemma~7.2]{Chen26}.
Related finite complex orbits appear in the work of Dubrovin and
Mazzocco \cite{DM00}.

For part~(3), let \(G_p=\mathrm{PSL}_2(\mathbb F_p)\).
An absolute \(G_p\)-structure is a generating pair of
\(G_p\) up to \(\operatorname{Aut}(G_p)\), which acts by
\(\mathrm{PGL}_2(\mathbb F_p)\)-conjugation
\cite[Lemma~2.1(a)]{Chen26}.
By \cite[Lemma~2.2]{Chen26}, every lift of such a pair
generates \(\mathrm{SL}_2(\mathbb F_p)\).
Changing the signs of the lifts induces the even sign
changes on trace triples and preserves the commutator.
Thus \cite[Lemma~2.1(b)]{Chen26} identifies \(S_D(p)/V\)
with the absolute \(G_p\)-structures of lifted commutator
trace \(D-2\). This identification is equivariant under
the Nielsen action.

By \cite[Lemma~2.1]{Chen26}, each
\(\operatorname{Out}^{+}(F_2)\)-orbit \(\Omega\) corresponds
to a geometric connected component of the moduli space
of absolute \(G_p\)-structures, whose degree over the
\(j\)-line equals \(|\Omega|\).
Chen's component-degree bound
\cite[Theorem~1.1]{Chen26} therefore yields
\(|\Omega|\ge p\).

We next recall the rotations used by Bourgain, Gamburd, and Sarnak
to convert orbit questions into multiplicative-order questions
\cite[\S~2.1]{BGS26}.
For $a\in\mathbb F_p\setminus\{\pm2\}$, write
$a=\lambda+\lambda^{-1}$ with $\lambda\in\mathbb F_{p^2}^{\times}$
and set $o(a)=\operatorname{ord}(\lambda)$.
This is independent of the choice between $\lambda$ and $\lambda^{-1}$.
We call $a$ \emph{split} or \emph{nonsplit}
according as $\chi_p(a^2-4)=1$ or $-1$; then $o(a)$ divides
$p-1$ or $p+1$, respectively.
The values $a=\pm2$ are called \emph{parabolic}.

The rotation $R_x(x,y,z)=(x,z,xz-y)\in H$ restricts on $x=a$
to $R_a(a,b,c)=(a,c,ac-b)$.
Diagonalizing its action on $(b,c)$ gives
\begin{equation}\label{eq:rotation-param}
\begin{aligned}
b(t)&=\alpha t+\beta t^{-1},\\
c(t)&=\lambda\alpha t+\lambda^{-1}\beta t^{-1},
\qquad t\in\langle\lambda\rangle,
\end{aligned}
\end{equation}
where
\begin{equation}\label{eq:kappa}
\alpha\beta=\kappa_a:=\frac{a^2-D}{a^2-4}.
\end{equation}
Indeed,
\[
b^2+c^2-abc=-(a^2-4)\alpha\beta=D-a^2.
\]

Since $D\ne4$, we have $\kappa_a\ne1$.
For a nonaxial point, $(\alpha,\beta)\ne(0,0)$.
If $R_a^j$ fixes the point, then $\lambda^j\alpha=\alpha$ and
$\lambda^{-j}\beta=\beta$, so $o(a)\mid j$.
Thus the rotation orbit has exactly $o(a)$ points.
In the split case $\alpha,\beta\in\mathbb F_p$; in the nonsplit case
$\beta=\alpha^p$ and $t$ lies in the norm-one torus.
Consequently, if $\kappa_a=0$, a nonaxial point is split, with
exactly one of $\alpha,\beta$ nonzero.

\section{Geometry and orbit structure}\label{sec:geometry}

\subsection{Finite-field curve estimates}\label{subsec:point-counts}

For $Q\in\{p-1,p+1\}$, let $T_Q$ be the split torus if $Q=p-1$
and the norm-one torus if $Q=p+1$.
Thus
\[
T_{p-1}(\mathbb F_p)=\mathbb F_p^\times,\qquad
T_{p+1}(\mathbb F_p)=
\{u\in\mathbb F_{p^2}^\times:u^{p+1}=1\}.
\]
The \emph{target torus} is $T_q$, with $q$ as in
\eqref{eq:choose-q}.

We make explicit the endgame of Bourgain, Gamburd, and Sarnak
\cite[\S\S~2--3]{BGS26}, using their final geometric arguments
\cite[\S\S~2.2 and 3.1--3.2]{BGS26} and keeping constants uniform in $D$.
Eddy et al.\ did this for $D=0$ \cite[Appendix]{EFLMT25}.
The first lemma counts intersections with torus subgroups;
the second finds a primitive fiber meeting two specified fibers.
All genus and boundary counts refer to smooth projective
normalizations.
We use the Hasse--Weil bound
\begin{equation}\label{eq:hasse-weil}
\bigl|\#C(\mathbb F_p)-(p+1)\bigr|\le2g(C)\sqrt p.
\end{equation}

\Needspace{10\baselineskip}
\begin{lemma}\label{lem:torus-points}
Let $p\ge7$, $Q_1,Q_2\in\{p-1,p+1\}$, and let $e\mid Q_1$
and $d\mid Q_2$.
Suppose $\alpha,\beta\in\mathbb F_{p^2}$,
$(\alpha,\beta)\ne(0,0)$, $\alpha\beta\ne1$, and
\[
\alpha,\beta\in\mathbb F_p\quad\text{if }Q_1=p-1,
\qquad
\beta=\alpha^p\quad\text{if }Q_1=p+1.
\]
Let $N$ be the number of pairs
$(u,v)\in T_{Q_1}(\mathbb F_p)\times T_{Q_2}(\mathbb F_p)$
satisfying
\begin{equation}\label{eq:torus-product}
\alpha u^e+\beta u^{-e}=v^d+v^{-d}=:y,
\qquad y^2\ne D.
\end{equation}
Then
\begin{equation}\label{eq:torus-count}
\left|\frac{N}{ed}-\frac{p}{ed}\right|
\le4\sqrt p+8.
\end{equation}
\end{lemma}

\begin{proof}
The stated reality conditions define an $\mathbb F_p$-curve in
$T_{Q_1}\times T_{Q_2}$.
Set $h=\gcd(e,d)$.
If $\alpha\beta\notin\{0,1\}$, the two covers of the $y$-line have
degrees $2e$ and $2d$, with disjoint finite branch loci
$\{\pm2\sqrt{\alpha\beta}\}$ and $\{\pm2\}$.
The tame fiber-product argument in the proof of
\cite[Lemma~12]{BGS26}, which requires only prime-to-$p$ exponents,
gives a geometrically irreducible cover of degree $4ed$.
Its four finite branch values have inertia order $2$, while above
infinity there are $4h$ points of inertia order $ed/h$.
Riemann--Hurwitz yields
\[
2g-2=-8ed+8ed+4(ed-h),\qquad g=2ed-2h+1.
\]
The affine fiber product is smooth; removing its points at infinity
and those with $y^2=D$ deletes at most $4h+8ed$ geometric points.
Hasse--Weil gives
\[
|N-p|\le(4ed-4h+2)\sqrt p+4h+8ed+1,
\]
which implies \eqref{eq:torus-count} for $p\ge7$.

If $\alpha\beta=0$, replacing $u$ by $u^{-1}$ if necessary gives
$y=\alpha u^e$ with $\alpha\ne0$.
Its finite branch locus is contained in $\{0\}$ and hence disjoint
from $\{\pm2\}$; the same argument gives degree $2ed$,
genus $ed-d-h+1$, and at most $2d+2h$ boundary points.
Deleting $y^2=D$ costs at most $4ed$ further points, giving the
stronger error $2\sqrt p+4$ after division by $ed$.
These geometric calculations also apply to the split and nonsplit
$\mathbb F_p$-forms.
\end{proof}

To connect two primitive fibers, we need a middle coordinate
that occurs on both fibers and is itself primitive.
As in Bourgain, Gamburd, and Sarnak's incidence argument
\cite[Lemma~9 and \S~3.2]{BGS26}, this leads to two simultaneous
quadratic conditions.

\begin{lemma}\label{lem:fiber-intersections}
Let $a,c$ be good primitive traces and put
\begin{equation}\label{eq:qa}
q_a(y)=(a^2-4)y^2+4(D-a^2).
\end{equation}
For $d\mid q$, let
$K_d=\{t^d:t\in T_q(\mathbb F_p)\}$, the subgroup of index $d$.
Let $G_d$ be the number of triples
$(s,v,w)\in K_d\times\mathbb F_p^2$ such that, with $y=s+s^{-1}$,
\begin{equation}\label{eq:intersection-system}
v^2=q_a(y),\qquad w^2=q_c(y),\qquad y^2\ne D.
\end{equation}
When $a^2=c^2$, impose only $v^2=q_a(y)$ and count pairs $(s,v)$.
Then
\begin{equation}\label{eq:fiber-count}
\left|G_d-\frac pd\right|\le16\sqrt p+16.
\end{equation}
\end{lemma}

\begin{proof}
Let $N_d$ be the number of solutions to
\eqref{eq:intersection-system} with $y=t^d+t^{-d}$ and
$t\in T_q(\mathbb F_p)$, omitting $w$ when $a^2=c^2$.
The power map $t\mapsto t^d$ is $d$-to-one onto $K_d$,
so $N_d=dG_d$.

Suppose first that $a^2\ne c^2$.
The branch pair of $v^2=q_a(y)$ is
$\{\pm2\sqrt{\kappa_a}\}$, nonzero and disjoint from $\{\pm2\}$.
Since $r\mapsto(r-D)/(r-4)$ is injective, the branch pairs for $a,c$
are also disjoint.
The argument of \cite[Lemma~9]{BGS26} gives a geometrically
irreducible degree-four cover defined by the two quadratic conditions.
After pulling back by $y=t^d+t^{-d}$, the two radicands have disjoint
simple zeros, so their square classes in
$\overline{\mathbb F}_p(t)$ remain independent.
Thus the full cover is geometrically irreducible of degree $8d$,
with six finite branch values of inertia order $2$ and eight points
above infinity of inertia order $d$.
Hence
\[
2g-2=-16d+24d+8(d-1),\qquad g=8d-3.
\]
After removing at most $8+16d$ points
at infinity and with $y^2=D$, Hasse--Weil gives
\[
\left|G_d-\frac pd\right|
\le\frac{(16d-6)\sqrt p+16d+9}{d}
\le16\sqrt p+16\qquad(p\ge7).
\]
If $a^2=c^2$, retain only one quadratic condition.
The cover then has degree $4d$, genus $2d-1$, and four points above
infinity, yielding the stronger error $4\sqrt p+8$.
Zeros of the quadratic discriminants are counted, and the same
argument applies to nonsplit forms.
\end{proof}

\subsection{Connectivity of primitive fibers}
\label{subsec:primitive-fibers}

The primitive-element sieve now gives an explicit uniform
Bourgain--Gamburd--Sarnak endgame \cite[\S~3]{BGS26,BGS26}.

\begin{proposition}[Primitive-fiber endgame]\label{prop:endgame}
Let $p\ge113$, $D\ne4$, and suppose
\begin{equation}\label{eq:endgame-criterion}
p>J(q)(16\sqrt p+16).
\end{equation}
Then $\mathcal C$ lies in a single $H$-orbit
$\mathcal U\subset S_D(p)$.

Moreover, $\mathcal U$ contains every nonaxial point with a coordinate
$\pm2$, and every nonaxial point having a coordinate of order $m>B_p$.

The hypothesis holds, in particular, for every prime $p>10^{12}$.
\end{proposition}

\begin{proof}
Consider a nonaxial rotation orbit of order $m$ in a torus of order
$Q\in\{p-1,p+1\}$, and put $e=Q/m$.
For each $d\mid q$, let $K_d$ be the index-$d$ subgroup of
$T_q(\mathbb F_p)$, and let $F_d$ count pairs
$(t,s)\in\langle\lambda\rangle\times K_d$ satisfying
\[
\alpha t+\beta t^{-1}=s+s^{-1}=:y,\qquad y^2\ne D.
\]
Since $u\mapsto u^e$ and $v\mapsto v^d$ have respectively $e$ and $d$
preimages on their image subgroups, the count $N$ in
Lemma~\ref{lem:torus-points} equals $edF_d$.
Thus
\[
\left|F_d-\frac{p}{ed}\right|\le4\sqrt p+8.
\]
Let $F_{\mathrm{prim}}$ denote the number of these pairs for which
$s$ has order $q$.
As $K_d$ is the subgroup of index $d$, \(s\in K_d\) precisely when
\(d\mid q/\operatorname{ord}(s)\), M\"obius inversion gives
\begin{equation}\label{eq:primitive-sieve}
\begin{aligned}
F_{\mathrm{prim}}
=\sum_{d\mid q}\mu(d)F_d\ge\frac{p}{e}\frac{\varphi(q)}q
-2^{\omega(q)}(4\sqrt p+8)>0
\end{aligned}
\end{equation}
whenever $m>B_p$.
Hence every such orbit reaches a good primitive fiber.

Each good primitive fiber is a single rotation orbit.
In the split case its points satisfy $\alpha\in\mathbb F_p^\times$
and $\beta=\kappa_a/\alpha$; in the nonsplit case they satisfy
$\alpha^{p+1}=\kappa_a$ and $\beta=\alpha^p$.
In either case multiplication by the primitive eigenvalue is
transitive on the $q$ points of the fiber.

For any two good primitive traces $a,c$, let $G_d$ be as in
Lemma~\ref{lem:fiber-intersections}, and let $G_{\mathrm{prim}}$
count the same configurations with $s$ of order $q$.
The same sieve yields
\[
G_{\mathrm{prim}}=\sum_{d\mid q}\mu(d)G_d
\ge\frac{p\varphi(q)}q
-2^{\omega(q)}(16\sqrt p+16)>0
\]
by \eqref{eq:endgame-criterion}.
The two square conditions give a point on each of the fibers with
fixed coordinates $a,c$, with the same good primitive middle coordinate.
These points are joined within that middle-coordinate fiber.
Hence all good primitive coordinate fibers belong to one $H$-orbit.

There is at least one good primitive trace:
$\varphi(q)\ge\sqrt{q/2}>4$ for $q\ge112$, whereas at most two of the
$\varphi(q)/2$ primitive traces satisfy $a^2=D$.
They are essential by Proposition~\ref{prop:literature}(2),
since $q>10$ and $a^2\ne D$ excludes axial points.

It remains to treat parabolic coordinates.
For $a=2$ the rotation acts as a nonzero translation on
$(b-c)^2=D-4$; for $a=-2$ its square acts as a nonzero translation on
$(b+c)^2=D-4$.
In either case a varying coordinate runs through $\mathbb F_p$
and hence attains a good primitive trace.
The final numerical assertion follows from
Lemma~\ref{lem:explicit-endgame}.
\end{proof}

\begin{lemma}
\label{lem:explicit-endgame}
For every integer $n\ge1$,
\begin{equation}\label{eq:J-power-bound}
J(n):=\frac{2^{\omega(n)}n}{\varphi(n)}<300n^{1/6}.
\end{equation}
Consequently, for every prime $p>10^{12}$ and
$q=p-\chi_p(-1)$,
\[
p>J(q)(16\sqrt p+16).
\]
\end{lemma}

\begin{proof}
The first bound is verified in Appendix~\ref{app:certificate}.
Since $q\le p+1<2p$ and $300\cdot2^{1/6}<400$,
we have $J(q)<400p^{1/6}$.
For $p>10^{12}$,
\[
J(q)(16\sqrt p+16)<6800p^{2/3}<p,
\]
as $p^{1/3}>10^4>6800$.
\end{proof}

\subsection{Sign invariance and descent to pure Vietas}
\label{subsec:sign-descent}

The choice \(4\mid q\) makes the good primitive traces
stable under negation. We use this symmetry to apply
the orbit-size bound of Proposition~\ref{prop:literature}(3)
to sign classes in the complement of \(\mathcal U\).

\begin{proof}[Proof of Theorem~\ref{thm:four-p}]
Proposition~\ref{prop:endgame} gives the single $H$-orbit
$\mathcal U=H\mathcal C\subseteq S_D(p)$.
If $\lambda$ has order $q$, then $-\lambda=\lambda^{q/2+1}$ also has
order $q$, since $4\mid q$ and $\gcd(q,q/2+1)=1$.
Negating a coordinate also preserves the condition $a^2\ne D$.
Thus $\mathcal C$ is $V$-invariant.

The Vieta involutions commute with $V$ because
$\epsilon_1\epsilon_2\epsilon_3=1$, and coordinate permutations
normalize $V$. Hence $H$ normalizes $V$, so
\[
V\mathcal U=V(H\mathcal C)=H(V\mathcal C)=H\mathcal C=\mathcal U.
\]
Since $S_D(p)$ and $\mathcal U$ are both $H$- and $V$-invariant,
so is $\mathcal B=S_D(p)\setminus\mathcal U$.

Suppose \(\mathcal B\ne\varnothing\).
Since \(\mathcal B\) is invariant under both \(H\) and \(V\),
the subset \(\mathcal B/V\subseteq S_D(p)/V\) is
invariant under the induced Nielsen action.
It therefore contains an entire
\(\operatorname{Out}^{+}(F_2)\)-orbit \(\Omega\), and
Proposition~\ref{prop:literature}(3) gives
\(|\Omega|\ge p\).

A nonidentity element of $V$ can fix a triple only if the two
coordinates whose signs it changes are zero. Such a triple is
nonessential by Proposition~\ref{prop:literature}(2), so $V$ acts
freely on $S_D(p)$. By $V$-invariance of $\mathcal B$, each class
in $\mathcal B/V$ therefore contains four points of $\mathcal B$.
Consequently,
\[
|\mathcal B|=4|\mathcal B/V|\ge4|\Omega|\ge4p.
\]
\end{proof}

We turn to the descent from $H$ to $\Gamma$.
Following Martin's normality and fixed-point strategy
\cite[proof of Theorem~4.7]{Martin25} and the odd-order construction
of Kingsbury-Neuschotz \cite[Lemmas~7.7--7.8]{KN26}, we find an
essential diagonal point whose coordinate transposition can be
realized using Vieta involutions.
The character-sum estimates below make the required existence
statement explicit for $p\ge113$.

\begin{proof}[Proof of Theorem~\ref{thm:descent}]
Write $\chi=\chi_p$, extended by $\chi(0)=0$, and put
$\Delta(a)=a^4-8a^2+4D$.
We seek an essential point $(x,a,a)$ with
\begin{equation}\label{eq:descent-character-conditions}
\chi(a+2)=1,\qquad
\chi(2-a)=-1,\qquad
\chi(\Delta(a))=1.
\end{equation}
Set
\[
W=\frac18\sum_{a\in\mathbb F_p}
(1+\chi(a+2))(1-\chi(2-a))(1+\chi(\Delta(a))).
\]
The odd terms cancel under $a\mapsto-a$, and
$\sum_a\chi(4-a^2)=-\chi(-1)$.
Hence
\begin{equation}\label{eq:W-expansion}
8W=p+\chi(-1)+\sum_a\chi(\Delta(a))
-\sum_a\chi((4-a^2)\Delta(a)).
\end{equation}
For $D\notin\{0,4\}$, the two character sums come from squarefree
polynomials of degrees $4$ and $6$; their smooth models have genera
$1$ and $2$, respectively.
Accounting for the points above infinity, Hasse--Weil yields
\[
\begin{aligned}
\sum_a\chi(\Delta(a))
&=-1+\varepsilon_1,
&|\varepsilon_1|&\le2\sqrt p,\\
\sum_a\chi((4-a^2)\Delta(a))
&=-\chi(-1)+\varepsilon_2,
&|\varepsilon_2|&\le4\sqrt p.
\end{aligned}
\]
For $D=0$, write $\Delta(a)=a^2(a^2-8)$.
The correction at $a=0$ and the genus-one estimate for
$(4-a^2)(a^2-8)$ give
\[
\begin{aligned}
\sum_a\chi(\Delta(a))
&=-1-\chi(-8),\\
\sum_a\chi((4-a^2)\Delta(a))
&=-\chi(-1)+\varepsilon-\chi(-32),
&|\varepsilon|&\le2\sqrt p.
\end{aligned}
\]
Since $\chi(-32)=\chi(-8)$, in both cases
\eqref{eq:W-expansion} implies
\begin{equation}\label{eq:W-general}
8W\ge p-3-6\sqrt p.
\end{equation}
The contribution from $a=-2$ is zero, and those from $a=2$ and the
at most four roots of $\Delta$ total at most $5/2$.
Thus the number of $a$ satisfying
\eqref{eq:descent-character-conditions} is at least
\[
W-\frac52\ge\frac{p-6\sqrt p-23}{8}>3\qquad(p\ge113).
\]
These candidates satisfy $a\ne0,\pm2$.
For each such $a$, the quadratic equation
\begin{equation}\label{eq:diagonal-quadratic}
x^2-a^2x+2a^2-D=0
\end{equation}
has two distinct roots.
Moreover, if $\rho^2=a+2$ and $t+t^{-1}=\rho$, then
$t^2+t^{-2}=a$ and $\chi(a-2)=-\chi(-1)$.
Thus $t$ lies in the split torus when $p\equiv3\pmod4$ and in the
norm-one torus when $p\equiv1\pmod4$.
In either case $o(a)=\operatorname{ord}(t^2)$ is odd, dividing
$(p-1)/2$ or $(p+1)/2$, respectively.
By Proposition~\ref{prop:literature}(2), a nonessential nonaxial point
with odd coordinate order has that order $3$ or $5$.
There are only $(\varphi(3)+\varphi(5))/2=3$ corresponding trace values.
Since there are more than three candidates $a$, one yields an
essential point $P=(x,a,a)$ with $o(a)$ odd.

Put $r=o(a)$.
On $z=a$, the rotation $R(x,y)=(y,ay-x)$ satisfies
$R^2=s_y s_x$, $R^r=1$, and $\tau_{xy}=s_yR$.
Thus the global element
\[
\gamma=s_y(s_y s_x)^{(r+1)/2}\in\Gamma
\]
satisfies $\gamma P=s_yRP=\tau_{xy}P=:Q=(a,x,a)$.
The transpositions $\tau_{yz}$ and $\tau_{xz}$ fix $P$ and $Q$,
respectively.
Since $\Gamma\triangleleft H$, both stabilize $\Gamma P$.
They generate $S_3$, so $\Gamma P$ is $H$-invariant and therefore
equals $S_D(p)$.
\end{proof}

\Needspace{11\baselineskip}
\section{Arithmetic bounds and the order obstruction}
\label{sec:arithmetic}

\subsection{Torus-subgroup intersections}
\label{subsec:torus-intersections}

Following Bourgain, Gamburd, and Sarnak
\cite[\S~4]{BGS26}, we use Corvaja and Zannier's gcd theorem
\cite[Theorem~2]{CZ13} to bound rotation-orbit intersections with
torus subgroups.
The following is the version used by Eddy et al.\
\cite[Lemma~2.1]{EFLMT25}, with the boundary and degrees recorded
for general $D$.

\begin{lemma}[Torus-subgroup intersection bound]
\label{lem:torus-intersection}
Let $\alpha,\beta\in\overline{\mathbb F}_p$ satisfy
$(\alpha,\beta)\ne(0,0)$ and $\alpha\beta\ne1$.
For positive integers $m,d$ with $p\nmid md$, write
$\mu_r=\{u\in\overline{\mathbb F}_p^\times:u^r=1\}$.
The number of solutions of
\begin{equation}\label{eq:torsion-intersection}
\alpha t+\beta t^{-1}=s+s^{-1},
\qquad (t,s)\in\mu_m\times\mu_d,
\end{equation}
is at most
\begin{equation}\label{eq:CZ-bound}
\max\left\{6\sqrt[3]{4}(md)^{1/3},\frac{48md}{p}\right\}.
\end{equation}
\end{lemma}

\begin{proof}
Let $N_{m,d}$ denote the number of solutions of
\eqref{eq:torsion-intersection}.
We follow \cite[proof of Lemma~2.1]{EFLMT25}, recording the boundary
explicitly.
If $\alpha\beta\notin\{0,1\}$, the projective closure is the smooth
genus-one cubic
\[
\alpha T^2S-TS^2-TZ^2+\beta SZ^2=0.
\]
Its boundary consists of
\[
[0:0:1],\quad [1:0:0],\quad [0:1:0],\quad [1:\alpha:0].
\]
The respective valuation pairs for $t=T/Z$ and $s=S/Z$ are
$(1,1)$, $(-1,1)$, $(1,-1)$, and $(-1,-1)$.
Thus $t,s$ are multiplicatively independent modulo constants.
Each has two simple poles, so $\deg t=\deg s=2$.
In the notation of \cite[Theorem~2]{CZ13},
\[
\chi=2g-2+4=4,\qquad
\deg(t^m)=2m,\quad \deg(s^d)=2d.
\]
The simple zeros and poles of $t,s$, together with $p\nmid md$,
ensure that $t^m,s^d$ have nonzero differentials.
Every solution contributes at least one to the gcd sum in that
theorem, so
\[
N_{m,d}\le
\max\left\{3(32md)^{1/3},\frac{48md}{p}\right\}.
\]
If $\alpha\beta=0$, symmetry lets us assume $\beta=0$, so
$t=(s+s^{-1})/\alpha$.
Its smooth projective model has genus zero, four boundary points,
$\deg t=2$, and $\deg s=1$.
The same argument gives the smaller bound
$\max\{3(8md)^{1/3},24md/p\}$.
\end{proof}

\subsection{The maximum-order obstruction}\label{subsec:obstruction}

Following Eddy et al.\ \cite[\S\S~2--3]{EFLMT25}, we choose a
coordinate of largest order outside $\mathcal U$.
The complement bound forces many coordinate values, while the
intersection estimate restricts their occurrence on a rotation orbit.

\Needspace{12\baselineskip}
\begin{proposition}\label{prop:obstruction}
Assume the hypotheses of Theorem~\ref{thm:four-p}, and let
$\mathcal B=S_D(p)\setminus\mathcal U\ne\varnothing$.
If $m$ is the largest eigenvalue order among all coordinates of
points of $\mathcal B$, then $m\mid p-1$ or $m\mid p+1$, and
\[
m\le B_p,\qquad
8p\le m^2K_p(m)^2.
\]
Moreover,
\[
m\le108K_p(m)^3
\qquad\text{or}\qquad
m\ge\frac{p}{24K_p(m)}.
\]
\end{proposition}

\begin{proof}
By Proposition~\ref{prop:endgame}, no point of $\mathcal B$ has a
parabolic coordinate and every nonparabolic coordinate of such a
point has eigenvalue order at most $B_p$.
Therefore $3\le m\le B_p$, and $m$ divides $p-1$ or $p+1$.
Put $K=K_p(m)$.

Every nonparabolic trace of order at most $m$ is associated with a
cyclic subgroup whose order is among
$\mathcal M_m(p-1)$ or $\mathcal M_m(p+1)$.
A subgroup of order $d$ contains fewer than $d/2$ nonparabolic trace
values, since inversion pairs its elements of order greater than two.
The split and nonsplit nonparabolic trace sets are disjoint.
Thus fewer than
\[
\frac12\sum_{\epsilon=\pm1}
\sum_{d\in\mathcal M_m(p+\epsilon)}d
\le\frac{mK}{2}
\]
trace values can occur in any coordinate of $\mathcal B$.
Any prescribed pair of coordinates determines at most two values of
the third.
Theorem~\ref{thm:four-p} therefore yields
\[
4p\le|\mathcal B|<2\left(\frac{mK}{2}\right)^2,
\qquad\text{so }\quad8p<m^2K^2.
\]

Choose a point of $\mathcal B$ with some coordinate of order $m$.
Its rotation orbit has exactly $m$ points and is contained in
$\mathcal B$.
Along this orbit another coordinate has the form
$\alpha t+\beta t^{-1}$ with $t\in\mu_m$.
Every such value is nonparabolic and has eigenvalue order at most $m$.
Hence its eigenvalues lie in a subgroup of order
$d\in\mathcal M_m(p-1)$ or $d\in\mathcal M_m(p+1)$.
Lemma~\ref{lem:torus-intersection} bounds the corresponding
number of pairs $(t,s)$ by
\[
\max\left\{6\sqrt[3]{4}(md)^{1/3},\frac{48md}{p}\right\}.
\]
Each nonparabolic trace has exactly two distinct inverse eigenvalues,
so the number of rotation parameters is at most half this bound.
Summing over the maximal divisors and using $d\le m$ gives
\[
m\le K\max\left\{3\sqrt[3]{4}m^{2/3},\frac{24m^2}p\right\}.
\]
If the first expression is maximal, then
$m^{1/3}\le3\sqrt[3]{4}K$, and therefore $m\le108K^3$.
Otherwise $m\le24Km^2/p$, so $m\ge p/(24K)$.
\end{proof}

\subsection{Weighted divisor antichains}\label{subsec:antichain}

To estimate $K_p(m)$ in Proposition~\ref{prop:obstruction},
we establish a weighted bound for divisor antichains.
The argument treats the exponents of $2$ and $3$ jointly and
exploits the coprimality of the parts of $p-1$ and $p+1$
supported on primes at least $5$.

\begin{proof}[Proof of Theorem~\ref{thm:antichain}]
Write $\sigma_{-s}(n)=\sum_{d\mid n}d^{-s}$.
A \emph{divisor antichain} is a set of divisors in which no two
distinct elements divide one another.
Define
\begin{equation}\label{eq:local-antichain-constants}
c_\ell=
\max_{a\ge0}
\ell^{-ta}\sum_{j=0}^{a}\ell^{-sj},
\qquad
C=\prod_{\ell\ge5\ {\rm prime}}c_\ell.
\end{equation}
Also put
\begin{equation}\label{eq:U-constant}
D_k=\sum_{j=0}^k2^{-sj}3^{-s(k-j)},
\qquad
U=\max\left\{1,\sup_{k\ge1}D_k6^{-tk}\right\}.
\end{equation}
Appendix~\ref{app:certificate} gives the finite verifications of
\begin{equation}\label{eq:certified-constants}
C<8\,600\,000,\qquad U<\frac94.
\end{equation}

Let $\mathcal A$ be any divisor antichain of an even integer
\[
n=2^a3^b r,\qquad (r,6)=1.
\]
Fix a divisor $v\mid r$.
The elements of $\mathcal A$ with remaining part $v$ correspond to
incomparable exponent pairs $(i,j)$.

If there are $k+1$ such pairs, list them with increasing $i$:
\[
(i_0,j_0),\ldots,(i_k,j_k),
\qquad i_0<i_1<\cdots<i_k.
\]
Their $j$-coordinates decrease strictly, so
\[
i_h\ge h,\qquad j_h\ge k-h.
\]
Therefore
\begin{equation}\label{eq:antichain-slice}
\sum_{h=0}^{k}2^{-si_h}3^{-sj_h}
\le
\sum_{h=0}^{k}2^{-sh}3^{-s(k-h)}
=D_k.
\end{equation}
For $k\ge1$, necessarily $a,b\ge k$, and hence
\[
2^{-ta}3^{-tb}D_k\le6^{-tk}D_k\le U.
\]
For $k=0$, the normalized contribution is at most $1\le U$.
This also handles $b=0$.
Summing the slices gives
\begin{equation}\label{eq:antichain-one-n}
\sum_{d\in\mathcal A}d^{-s}
\le U2^{ta}3^{tb}\sum_{v\mid r}v^{-s}
=U n^t\,\frac{\sigma_{-s}(r)}{r^t}.
\end{equation}

Now apply this estimate to the maximal-divisor antichains of
$P-1$ and $P+1$.
Write
\[
P\pm1=2^{a_\pm}3^{b_\pm}r_\pm,
\qquad (r_\pm,6)=1.
\]
Since $\gcd(P-1,P+1)=2$, we have $\gcd(r_-,r_+)=1$.

Set
\[
X=\frac{\sigma_{-s}(r_-)}{r_-^t},
\qquad
Y=\frac{\sigma_{-s}(r_+)}{r_+^t}.
\]
By multiplicativity and the definition of $C$,
\[
X\le C,\qquad Y\le C,
\]
and, since $r_-$ and $r_+$ are coprime,
\[
XY
=\frac{\sigma_{-s}(r_-r_+)}{(r_-r_+)^t}
\le C.
\]
It follows that
\[
X+Y\le C+1.
\]
Indeed, this is immediate if either $X\le1$ or $Y\le1$.
Otherwise,
\[
X+Y\le XY+1\le C+1.
\]

Applying the preceding antichain estimate to both $P-1$ and $P+1$,
and using $t>0$, gives
\begin{equation}\label{eq:antichain-pair}
\sum_{\epsilon=\pm1}
\ \sum_{d\in\mathcal M_m(P+\epsilon)}d^{-s}
\le U(P+1)^t(C+1).
\end{equation}

Finally, $d\le m$ implies
\[
d^\eta\le m^{\eta+s}d^{-s}.
\]
Using \eqref{eq:certified-constants},
\[
U(C+1)<\frac94(8\,600\,001)<19\,400\,000,
\]
which proves \eqref{eq:weighted-antichain}.
\end{proof}

\section{Proof of uniform transitivity}\label{sec:proof}

We finish by combining Chen's lower bound with the maximum-order
obstruction and the weighted maximal-divisor estimate.
The two order alternatives each yield a numerical contradiction,
which excludes a nonempty complement of $\mathcal U$.

\begin{proof}[Proof of Theorem~\ref{thm:main}]
Let $p>10^{320}$; throughout, choose $q=p-\chi_p(-1)$,
as in Section~\ref{sec:geometry}.

Suppose, for contradiction, that
$\mathcal B=S_D(p)\setminus\mathcal U$ is nonempty.
Let $m$ be its maximum coordinate order, and write $K=K_p(m)$.
Proposition~\ref{prop:obstruction} gives
\begin{equation}\label{eq:size-obstruction}
8p\le m^2K^2,
\end{equation}
and either
\begin{equation}\label{eq:small-order-case}
m\le108K^3
\end{equation}
or
\begin{equation}\label{eq:large-order-case}
p\le24mK.
\end{equation}

For the small-order alternative, by Theorem~\ref{thm:antichain},
\[
K<K_0(p+1)^t m^s.
\]
Combining this with \eqref{eq:small-order-case},
\begin{equation}\label{eq:small-m-bound}
m^{1-3s}
\le108K_0^3(p+1)^{3t}.
\end{equation}
On the other hand, \eqref{eq:size-obstruction} gives
\begin{equation}\label{eq:small-size-bound}
8p\le K_0^2(p+1)^{2t}m^{2+2s}.
\end{equation}
Raise \eqref{eq:small-size-bound} to $1-3s$ and eliminate $m$
using \eqref{eq:small-m-bound}.
The result is
\begin{equation}\label{eq:small-p-bound}
p^{97}
\le
108^{1090}K_0^{4000}8^{-365}
\left(1+\frac1p\right)^{268}.
\end{equation}

For $p>10^{320}$, we have $1+1/p<1001/1000$.
The exact comparison
\begin{equation}\label{eq:small-numerical-check}
10^{31040}8^{365}1000^{268}
>
108^{1090}K_0^{4000}1001^{268}
\end{equation}
contradicts \eqref{eq:small-p-bound}.

For the large-order alternative,
Lemma~\ref{lem:explicit-endgame} gives $J(q)<400p^{1/6}$.
Since $(p+1)/p<2$ and $4\sqrt p+8<5\sqrt p$, we get
\begin{equation}\label{eq:large-B}
B_p<4000p^{2/3}.
\end{equation}

Using \eqref{eq:K-bound}, $p+1<2p$ and $K_0<2\cdot10^7$,
we obtain
\[
K<4\cdot10^7p^tm^s,\qquad m<4000p^{2/3}.
\]
Therefore \eqref{eq:large-order-case} implies
\[
\begin{aligned}
p<24(4\cdot10^7)(4000)^{1+s}
  p^{t+\frac23(1+s)}<10^{13}p^{2381/3000}.
\end{aligned}
\]
Hence
\[
p^{619/3000}<10^{13},
\]
which is impossible for $p>10^{320}$.

Both alternatives are excluded.
Thus $H$ is transitive on $S_D(p)$.
Theorem~\ref{thm:descent} gives $\Gamma$-transitivity.
\end{proof}

\appendix

\section{Finite numerical verifications}\label{app:certificate}

We justify the bounds used above:
\begin{equation}\label{eq:appendix-goals}
C<8\,600\,000,\qquad U<\frac94,\qquad
J(n)<300n^{1/6}\quad(n\ge1),
\end{equation}
where $C,U$ are defined in
\eqref{eq:local-antichain-constants}--\eqref{eq:U-constant}.

Put $f_\ell(a)=\ell^{-ta}\sum_{j=0}^a\ell^{-sj}$.
Since the summands decrease with $j$,
\[
\frac{f_\ell(a+1)}{f_\ell(a)}
\le\ell^{-t}\frac{a+2}{a+1}.
\]
The inequalities
\[
23^{1000}<22^{1000}2^{67},\qquad
1024^{-t}(1+1024^{-s})<1
\]
show that $f_\ell$ decreases for $a\ge21$, and that
$c_\ell=1$ for all primes $\ell\ge1024$.
Thus $C$ is determined by primes $5\le\ell<1024$ and
$0\le a\le21$.

For $U$, let $R=2^{-s}6^{-t}$, $S=3^{-s}6^{-t}$ and
$h_k=\sum_{j=0}^kR^jS^{k-j}$.
Since $R>S$ and $R<9/10$,
\[
h_{k+1}=Rh_k+S^{k+1}
\le R\frac{k+2}{k+1}h_k<h_k\qquad(k\ge9).
\]
Hence $U=\max\{1,h_1,\dots,h_9\}$.

For $n=\prod_\ell\ell^{a_\ell}$,
\[
\frac{J(n)^6}{n}
=\prod_{\ell\mid n}
\frac{(2\ell)^6}{\ell^{a_\ell}(\ell-1)^6}
\le\prod_{\ell\mid n}
\frac{(2\ell)^6}{\ell(\ell-1)^6}.
\]
The function $2x^{5/6}/(x-1)$ decreases for $x>1$.
The local factors are greater than $1$ for $\ell\le67$ and less
than $1$ for $\ell\ge71$.
Extending the remaining product to all primes at most $67$ gives
\begin{equation}\label{eq:J-finite-product}
\prod_{\substack{\ell\le67\\\ell\ \mathrm{prime}}}
\frac{(2\ell)^6}{\ell(\ell-1)^6}
<287^6<300^6,
\end{equation}
proving \eqref{eq:J-power-bound} by a finite integer comparison.

The remaining finite comparisons yield the rational upper bounds
\[
\begin{aligned}
C<8\,600\,000,\qquad U<\frac94.
\end{aligned}
\]

\section*{Acknowledgments}

The author is grateful to Daniel E. Martin for his guidance,
helpful discussions, and valuable comments on earlier
versions of this manuscript.

\medskip
\noindent\textbf{AI assistance.}
OpenAI's ChatGPT was used to assist with numerical computations, literature searches, identification of references, and the writing and
revision of the manuscript. The author takes full
responsibility for the mathematical arguments and results.

\end{document}